\documentclass[review]{elsarticle}
\usepackage{lineno,hyperref}
\usepackage{amsmath}
\usepackage{multirow}
\usepackage{mdwlist}
\usepackage{xcolor}
\usepackage{subfigure}
\usepackage{colortbl}
\usepackage{amssymb}
\modulolinenumbers[5]

\journal{Journal of Computational and Applied Mathematics}

\biboptions{numbers,sort&compress}

\newtheorem{Theorem}{Theorem}
\newtheorem{Example}{Example}
\newtheorem{Lemma}{Lemma}
\newtheorem{Remark}{Remark}
\newtheorem{Corollary}{Corollary}
\newtheorem{Definition}{Definition}
\newtheorem{Algorithm}{Algorithm}

\newenvironment{Proof}[1][Proof.]{\begin{trivlist}
\item[\hskip \labelsep {\bfseries #1}]}{\ep\end{trivlist}}
\newcommand{\ep}{\hfill\rule{0.15cm}{0.35cm}\vskip 0.3cm}

\begin{document}

\begin{frontmatter}

\title{Fully Finite Element Adaptive Algebraic Multigrid Method for Time-Space Caputo-Riesz Fractional Diffusion Equations}

\author[mymainaddress]{Xiaoqiang Yue}

\author[mymainaddress]{Weiping Bu}

\author[mysecondaryaddress]{Shi Shu\corref{mycorrespondingauthor}}
\cortext[mycorrespondingauthor]{Corresponding author}
\ead{shushi@xtu.edu.cn}

\author[mymainaddress]{Menghuan Liu}

\author[mymainaddress]{Shuai Wang}

\address[mymainaddress]{School of Mathematics and Computational Science, Xiangtan University, Hunan 411105, P.R. China}
\address[mysecondaryaddress]{Hunan Key Laboratory for Computation and Simulation in Science and Engineering, Xiangtan University, Hunan 411105, P.R. China}

\begin{abstract}
The paper aims to establish a fully discrete finite element (FE) scheme and provide cost-effective solutions for 
one-dimensional time-space Caputo-Riesz fractional diffusion equations on a bounded domain $\Omega$.
Firstly, we construct a fully discrete scheme of the linear FE method in both temporal and spatial directions,
derive many characterizations on the coefficient matrix
and numerically verify that the fully FE approximation possesses the saturation error order under $L^2(\Omega)$ norm.
Secondly, we theoretically prove the estimation $1+\mathcal{O}(\tau^\alpha h^{-2\beta})$ on the condition number of the coefficient matrix,
in which $\tau$ and $h$ respectively denote time and space step sizes. Finally, on the grounds of the estimation and fast Fourier transform,
we develop and analyze an adaptive algebraic multigrid (AMG) method with low algorithmic complexity,
reveal a reference formula to measure the strength-of-connection tolerance
which severely affect the robustness of AMG methods in handling fractional diffusion equations,
and illustrate the well robustness and high efficiency of the proposed algorithm compared with the classical AMG,
conjugate gradient and Jacobi iterative methods.
\end{abstract}


\end{frontmatter}

\linenumbers

\section{Introduction}

In recent years, there has been an explosion of research interest in numerical solutions for fractional differential equations,
mainly due to the following two aspects: (i) the huge majority can't be solved analytically,
(ii) the analytical solution (if luckily derived)
always involve certain infinite series which sharply drives up the costs of its evaluation.
Various numerical methods have been proposed to approximate
more accurately and faster, such as finite difference (FD) method \cite{l-002,d-001,g-001,y-007,c-002,w-002,c-001,w-003,w-001,y-004},
finite element (FE) method \cite{e-001,z-002,b-002,b-001,m-001,f-001,b-004},
finite volume \cite{l-003} method and spectral (element) method \cite{l-001,z-001,y-005,y-006,y-001,y-002,y-003}.
An essential challenge against standard differential equations
lies in the presence of the fractional differential operator,
which gives rise to nonlocality (space fractional, nearly dense or full coefficient matrix) or
memory-requirement (time fractional, the entire time history of evaluations) issue, resulting in a vast computational cost.

Preconditioned Krylov subspace methods are regarded as one of the potential solutions to the aforementioned challenge.
Numerous preconditioners with various Krylov-subspace methods have been constructed respectively for one- and two-dimensional,
linear and nonlinear space-fractional diffusion equations (SFDE) \cite{l-004,m-002,j-001,g-002,d-002}.
Multigrid method has been proven to be a superior solver
and preconditioner for ill-conditioned Toeplitz systems as well as SFDE.
Pang and Sun propose an efficient and robust geometric multigrid (GMG)
with fast Fourier transform (FFT) for one-dimensional SFDE by an implicit FD scheme \cite{p-001}.
Bu et al. employ the GMG to one-dimensional multi-term time-fractional advection-diffusion equations
via a fully discrete scheme by FD method in temporal and FE method in spatial directions \cite{b-003}.
Jiang and Xu construct optimal GMG for two-dimensional SFDE to get FE approximations \cite{j-002}.
Chen et al. make the first attempt to present an algebraic multigrid (AMG) method with line smoothers to
the fractional Laplacian through localizing it into a nonuniform elliptic equation \cite{c-003}.
Zhao et al. invoke GMG for one-dimensional Riesz SFDE by an adaptive FE scheme using hierarchical matrices \cite{z-003}.
From the survey of references, in spite of quite a number of contributions to numerical methods and preconditioners,
there are no calculations taking into account of fully discrete FE schemes and AMG methods
for time-space Caputo-Riesz fractional diffusion equations.

In this paper, we are concerned with the following time-space Caputo-Riesz fractional diffusion equation (CR-FDE)
\begin{align}\label{chp-01-01}
  &{}_0^CD^\alpha_tu(x,t)=\frac{\partial^{2\beta}u(x,t)}{\partial|x|^{2\beta}}+f(x,t),~~
  t\in I=(0,T],~x\in\Omega=(a,b) \\
  &\label{chp-01-02}u(x,t)=0,~~t\in I,~x\in\partial\Omega \\
  &\label{chp-01-03}u(x,0)=\psi_0(x),~~x\in\Omega
\end{align}
with orders $\alpha\in(0,1)$ and $\beta\in(1/2,1)$, the Caputo and Riesz fractional derivatives are respectively defined by
\begin{eqnarray*}
  {}_0^CD^\alpha_tu=\frac{1}{\Gamma(1-\alpha)}\int_0^t(t-s)^{-\alpha}\frac{\partial u}{\partial s}ds,~
  \frac{\partial^{2\beta}u}{\partial|x|^{2\beta}}=-\frac{1}{2\cos(\beta\pi)}({}_xD_L^{2\beta}u
  + {}_xD_R^{2\beta}u),
\end{eqnarray*}
where
\begin{eqnarray*}
  {}_xD_L^{2\beta}u = \frac{1}{\Gamma(2-2\beta)}\frac{\partial^2}{\partial x^2}\int_a^x(x-s)^{1-2\beta}uds,~
  {}_xD_R^{2\beta}u = \frac{1}{\Gamma(2-2\beta)}\frac{\partial^2}{\partial x^2}\int_x^b(s-x)^{1-2\beta}uds.
\end{eqnarray*}

The remainder of this paper proceeds as follows. A fully discrete FE method of \eqref{chp-01-01}-\eqref{chp-01-03} is developed in Section 2.
Section 3 comes up with the theoretical estimation and verification experiments on the condition number of the coefficient matrix.
The classical AMG method is introduced in Section 4 followed by its uniform convergence analysis
and the construction of an adaptive AMG method. Section 5 reports and analyzes numerical results to show the benefits.
We close in Section 6 with some concluding remarks.

\section{Fully discrete finite element scheme for the CR-FDE}

For simplicity, following \cite{x-001}, we will use the symbols $\lesssim$, $\gtrsim$ and $\simeq$ throughout the paper.
$u_1\lesssim v_1$ means $u_1\leq C_1v_1$, $u_2\gtrsim v_2$ means $u_2\ge c_2v_2$ while $u_3\simeq v_3$ means
$c_3v_3\leq u_3\leq C_3v_3$, where $C_1$, $c_2$, $c_3$ and $C_3$ are generic positive constants independent of variables, time and space step sizes.

\subsection{Reminder about fractional calculus}

In this subsection, we briefly introduce some fractional derivative spaces and several auxiliary results.
Here the $L^2$ inner product and norm are denoted by
$$(u,v)_{L^2(\Omega)}=\int_{\Omega}uvdx,~\|u\|_{L^2(\Omega)}=(u,u)^{\frac{1}{2}}_{L^2(\Omega)}.$$

\begin{Definition} (Left and right fractional derivative spaces)
  For constant $\mu>0$, define norms
\begin{eqnarray*}
  \|u\|_{J_L^\mu(\Omega)} := (\|u\|^2_{L^2(\Omega)}+\|{}_xD_L^\mu u\|^2_{L^2(\Omega)})^{\frac{1}{2}},~
  \|u\|_{J_R^\mu(\Omega)} := (\|u\|^2_{L^2(\Omega)}+\|{}_xD_R^\mu u\|^2_{L^2(\Omega)})^{\frac{1}{2}},
\end{eqnarray*}
  and let $J_{L,0}^\mu(\Omega)$ and $J_{R,0}^\mu(\Omega)$ be closures of $C_0^\infty(\Omega)$
  under $\|\cdot\|_{J_L^\mu(\Omega)}$ and $\|\cdot\|_{J_R^\mu(\Omega)}$, respectively.
\end{Definition}

\begin{Definition}\label{def-01} (Fractional Sobolev space)
  For constant $\mu>0$, define the norm
\begin{eqnarray}\label{chp-03-005}
  \|u\|_{H^\mu(\Omega)} := (\|u\|^2_{L^2(\Omega)}+\||\xi|^\mu\tilde{u}\|^2_{L^2(\Omega_\xi)})^{\frac{1}{2}},
\end{eqnarray}
  and let $H_0^\mu(\Omega)$ be the closure of $C_0^\infty(\Omega)$ under $\|\cdot\|_{H^\mu(\Omega)}$,
  where $\tilde{u}$ is the Fourier transform of $u$.
\end{Definition}

\begin{Remark}
  The equivalence between \eqref{chp-03-005} and the general definition of the norm has been established in \cite{b-001},
  which implies the reasonability of Definition \ref{def-01}.
\end{Remark}

\begin{Lemma} (see \cite{z-002}, Proposition 1)\label{clp-02-01}
  If constant $\mu\in(0,1)$, $u,v\in J_{L,0}^{2\mu}(\Omega)$ (or $J_{R,0}^{2\mu}(\Omega)$), then
\begin{eqnarray*}
  ({}_xD_L^{2\mu}u,v)_{L^2(\Omega)}=({}_xD_L^\mu u,{}_xD_R^\mu v)_{L^2(\Omega)},~
  ({}_xD_R^{2\mu}u,v)_{L^2(\Omega)}=({}_xD_R^\mu u,{}_xD_L^\mu v)_{L^2(\Omega)}.
\end{eqnarray*}
\end{Lemma}

\begin{Lemma} (see \cite{e-001}, Lemma 2.4)\label{clp-02-02}
  For constant $\mu > 0$, we have
\begin{eqnarray}\label{chp-03-002}
  ({}_xD_L^\mu u, {}_xD_R^\mu u)_{L^2(\Omega)} = \cos(\pi \mu)\|{}_xD_L^\mu u\|^2_{L^2(\Omega)}.
\end{eqnarray}
\end{Lemma}

\begin{Lemma} (Fractional Poincar{\'e}-Friedrichs inequality, see \cite{e-001}, Theorem 2.10)\label{clp-02-03}
  For $u\in J_{L,0}^\mu(\Omega)$, we have
\begin{eqnarray}\label{chp-03-003}
  \|u\|_{L^2(\Omega)}\lesssim \|{}_xD_L^\mu u\|_{L^2(\Omega)}.
\end{eqnarray}
\end{Lemma}

\subsection{Derivation of the fully discrete scheme}

By Lemma \ref{clp-02-01}, we get the variational (weak) formulation of \eqref{chp-01-01}-\eqref{chp-01-03}:
given $f\in L^2(\Omega,I)$, $\phi_0\in L^2(\Omega)$ and $Q_t:=\Omega\times(0,t)$, to find $u\in \mathcal{H}$ subject to $u(x,0)=\psi_0(x)$ and
\begin{eqnarray}\label{chp-02-003}
  \Big{(}{}_0^CD^\alpha_\sigma u,v\Big{)}_{Q_t}+B^t_\Omega(u,v)=(f,v)_{Q_t},~\forall v\in \mathcal{H}^*,
\end{eqnarray}
where $\mathcal{H}:=H_0^\beta(\Omega)\times H^1(I)$, $\mathcal{H}^*:=H_0^\beta(\Omega)\times L^2(I)$, and
\begin{align*}
  &\Big{(}{}_0^CD^\alpha_\sigma u,v\Big{)}_{Q_t}=\int_0^t\Big{(}{}_0^CD^\alpha_\sigma u,v\Big{)}_{L^2(\Omega)}d\sigma,
  ~\Big{(}f,v\Big{)}_{Q_t}=\int_0^t\Big{(}f,v\Big{)}_{L^2(\Omega)}d\sigma, \notag \\
  &B^t_\Omega(u,v)=\int_0^t\frac{1}{2\cos(\beta\pi)}\Big{[}({}_xD^\beta_Lu,{}_xD^\beta_Rv)_{L^2(\Omega)}
  +({}_xD^\beta_Ru,{}_xD^\beta_Lv)_{L^2(\Omega)}\Big{]}d\sigma. \notag
\end{align*}


In order to acquire numerical solutions of $u$, we firstly make a (possibly nonuniform) temporal discretization by points $0=t_0<t_1<\cdots<t_N=T$,
and a uniform spatial discretization by points $x_i = a+ih$ ($i=0,1,\cdots,M$), where $h=(b-a)/M$ represents the space step size. Let
\begin{eqnarray*}
I_j=(t_{j-1},t_j),~\tilde{I}_j=(0,t_j),~j=1,2,\cdots,N;~
\Omega_h = \{\Omega_l:\Omega_l=(x_{l-1},x_l),~l=1,2,\cdots,M\}.
\end{eqnarray*}

We observe that it is convenient to form the FE spaces in tensor products
\begin{eqnarray*}
  \mathcal{V}_n=\mathcal{V}^\beta_h(\Omega_h)\times\mathcal{V}_\tau(\tilde{I}_n),~
  \mathcal{V}^*_n=\mathcal{V}^\beta_h(\Omega_h)\times\mathcal{V}^*_\tau(I_n),
\end{eqnarray*}
where
\begin{eqnarray*}
  & \mathcal{V}^\beta_h(\Omega_h)=\{w_h\in H_0^\beta(\Omega)\cap C(\bar{\Omega}):~
  w_h(x)\big{|}_{\Omega_l}\in\mathcal{P}_1(\Omega_l),~l=1,\cdots,M\},\\
  & \mathcal{V}_\tau(\tilde{I}_n) = \{v_\tau\in \mathcal{C}(\overline{\tilde{I}_n}):~v_\tau(0)=1,~
  v_\tau(t)\big{|}_{I_j}\in\mathcal{P}_1(I_j),~j=1,\cdots, n\},\\
  & \mathcal{V}^*_\tau(I_n) = \{v_\tau\in L^2(I_n):~
  v_\tau(t)\big{|}_{I_n}\in\mathcal{P}_0(I_n)\},
\end{eqnarray*}
and $\mathcal{P}_k$ denotes the set of all polynomials of degree $\leq k$.

\begin{Remark}
  Apparently, for a given $u_{h\tau}(x,t)\in\mathcal{V}_n$, we have $\partial u_{h\tau}/\partial t\in \mathcal{V}^*_n$,
  where $\partial u_{h\tau}/\partial t$ is obtained by differentiating $u_{h\tau}$
  with respect to $t$ on each subinterval $I_j$ ($j=1,2,\cdots,N$).
\end{Remark}

We obtain a fully discrete FE scheme in temporal and spatial directions of problem \eqref{chp-02-003}:
given $Q_n:=\Omega_h\times I_n$, to find $u_{h\tau}\in \mathcal{V}_n$ such that $u_{h\tau}(x,0)=\psi_{0,I}(x)$ and
\begin{eqnarray}\label{chp-02-004}
  \Big{(}{}_0^CD^\alpha_t u_{h\tau},v_{h\tau}\Big{)}_{Q_n}+B^n_\Omega(u_{h\tau},v_{h\tau})=
  (f,v_{h\tau})_{Q_n},~\forall v_{h\tau}\in\mathcal{V}^*_n,
\end{eqnarray}
where $\psi_{0,I}(x)\in\mathcal{V}_n$ satisfying $\psi_{0,I}(x_i)=\psi_0(x_i)$ ($i=0,1,\cdots,M$), and
\begin{align*}
  & \Big{(}{}_0^CD^\alpha_tu_{h\tau},v_{h\tau}\Big{)}_{Q_n}=\int_{t_{n-1}}^{t_n} ({}_0^CD^\alpha_t u_{h\tau}, v_{h\tau})_{L^2(\Omega)}dt,~
  (f,v_{h\tau})_{Q_n}=\int_{t_{n-1}}^{t_n}(f,v_{h\tau})_{L^2(\Omega)}dt,\\
  & B^n_\Omega(u_{h\tau},v_{h\tau})
  = \int_{t_{n-1}}^{t_n}\frac{1}{2\cos(\beta\pi)}\Big{[}({}_xD^\beta_Lu_{h\tau},{}_xD^\beta_Rv_{h\tau})_{L^2(\Omega)}
  +({}_xD^\beta_Ru_{h\tau},{}_xD^\beta_Lv_{h\tau})_{L^2(\Omega)}\Big{]}dt.
\end{align*}

Let
\begin{eqnarray*}
  \mathcal{L}_0(t)=\left \{
    \begin{aligned}
    & \frac{t_1-t}{\tau_1},~t \in I_1 \\
    & 0,\qquad ~t\in \tilde{I}_n\setminus I_1
    \end{aligned}
  \right.,~
  \tilde{\mathcal{L}}_0(t)=\frac{1}{\Gamma(1-\alpha)}\int_{t_0}^t\frac{d\mathcal{L}_0(s)}{(t-s)^\alpha},~
  \hat{\mathcal{L}}_0(t)=\frac{1}{\Gamma(1-\alpha)}\int_{t_0}^{t_1}\frac{d\mathcal{L}_0(s)}{(t-s)^\alpha},
\end{eqnarray*}
\begin{eqnarray*}
    \mathcal{L}_k(t)=\left \{
    \begin{aligned}
    & \frac{t_{k+1}-t}{\tau_{k+1}},~t \in I_{k+1} \\
    & \frac{t-t_{k-1}}{\tau_k},~t \in I_k \\
    & 0,\qquad \quad ~~t\in \tilde{I}_n\setminus (I_k \cup I_{k+1})
    \end{aligned}
  \right.,~
  \hat{\mathcal{L}}_k(t)=\frac{1}{\Gamma(1-\alpha)}\int_{t_{k-1}}^{t_{k+1}}\frac{d\mathcal{L}_k(s)}{(t-s)^\alpha},~
  k=1,\cdots,n-1
\end{eqnarray*}
and
\begin{eqnarray*}
  \mathcal{L}_n(t)=\left \{
    \begin{aligned}
    & \frac{t-t_{n-1}}{\tau_n},~t \in I_n \\
    & 0,\qquad \quad ~~t\in \tilde{I}_n\setminus I_n
    \end{aligned}
  \right.,~\tilde{\mathcal{L}}_n(t)=\frac{1}{\Gamma(1-\alpha)}\int_{t_{n-1}}^t\frac{d\mathcal{L}_n(s)}{(t-s)^\alpha}.
\end{eqnarray*}
Note that
\begin{eqnarray*}
  \mathcal{V}^*_n=\textsf{span}\{\phi_l(x)\times 1,~l=1,\cdots,M-1\},
\end{eqnarray*}
where $\phi_l(x)$ is the shape function at $x_l\in\Omega_h$. Using
\begin{eqnarray*}\label{chp-02-005}
  u_{h\tau}(x,t)=u_h^0(x)\mathcal{L}_0(t)+\sum_{k=1}^{n-1}u_h^k(x)\mathcal{L}_k(t)+u_h^n(x)\mathcal{L}_n(t),
\end{eqnarray*}
we have
\begin{align}\label{chp-03-006}
  &({}_0^CD^\alpha_tu_{h\tau},\phi_l\times 1)_{Q_1}
  =(u_h^0,\phi_l)_{L^2(\Omega)}(\tilde{\mathcal{L}}_0,1)_{L^2(I_1)}
  +(u_h^1,\phi_l)_{L^2(\Omega)}(\tilde{\mathcal{L}}_1,1)_{L^2(I_1)}, \\ \label{chp-03-007}
  &({}_0^CD^\alpha_tu_{h\tau},\phi_l\times 1)_{Q_n}
  =(u_h^0,\phi_l)_{L^2(\Omega)}(\hat{\mathcal{L}}_0,1)_{L^2(I_n)}+ \notag \\
  &\qquad \sum_{k=1}^{n-1}(u_h^k,\phi_l)_{L^2(\Omega)}(\hat{\mathcal{L}}_k,1)_{L^2(I_n)}
  +(u_h^n,\phi_l)_{L^2(\Omega)}(\tilde{\mathcal{L}}_n,1)_{L^2(I_n)},~n>1, \\ \label{chp-03-014}
  & \int_{t_{n-1}}^{t_n}1\times({}_xD^{\beta}_Lu_{h\tau},{}_xD^{\beta}_R\phi_l)_{L^2(\Omega)}dt
  =({}_xD^{\beta}_Lu_h^0,{}_xD^{\beta}_R\phi_l)_{L^2(\Omega)}
  (\mathcal{L}_0,1)_{L^2(I_n)}+\notag\\
  &\qquad \sum_{k=1}^{n-1}({}_xD^{\beta}_Lu_h^k,{}_xD^{\beta}_R\phi_l)_{L^2(\Omega)}
  (\mathcal{L}_k,1)_{L^2(I_n)} +
  ({}_xD^{\beta}_Lu_h^n,{}_xD^{\beta}_R\phi_l)_{L^2(\Omega)}
  (\mathcal{L}_n,1)_{L^2(I_n)}, \\ \label{chp-03-008}
  &\int_{t_{n-1}}^{t_n}1\times
  ({}_xD^{\beta}_Ru_{h\tau},{}_xD^{\beta}_L\phi_l)_{L^2(\Omega)}dt
  =({}_xD^{\beta}_Ru_h^0,{}_xD^{\beta}_L\phi_l)_{L^2(\Omega)}
  (\mathcal{L}_0,1)_{L^2(I_n)}+ \notag\\
  &\qquad \sum_{k=1}^{n-1}({}_xD^{\beta}_Ru_h^k,{}_xD^{\beta}_L\phi_l)_{L^2(\Omega)}
  (\mathcal{L}_k,1)_{L^2(I_n)}+
  ({}_xD^{\beta}_Ru_h^n,{}_xD^{\beta}_L\phi_l)_{L^2(\Omega)}
  (\mathcal{L}_n,1)_{L^2(I_n)}.
\end{align}

Substituting \eqref{chp-03-006}-\eqref{chp-03-008} into \eqref{chp-02-004}, yields
\begin{eqnarray}\label{chp-02-009}
  C^n_{h\tau}U^n_{h\tau} = G^n_{h\tau},
\end{eqnarray}
where the coefficient matrix
\begin{eqnarray}\label{chp-02-008}
  C^n_{h\tau} = M_h+\frac{\Gamma(3-\alpha)}{2}\tau^\alpha_n A_h^\beta,
\end{eqnarray}
the right-hand side vector
\begin{align*}
  &G^n_{h\tau}=\Gamma(3-\alpha)\tau_n^{\alpha-1}F^n_{h\tau}+
  \Big{[}M_h-\frac{\Gamma(3-\alpha)}{2}\tau^\alpha_n A_h^\beta\Big{]}U^{n-1}_{h\tau} - \sum_{k=1}^{n-1} \tau_n^{\alpha-1} \times \notag \\
  &\quad \frac{(t_n-t_{k-1})^{2-\alpha}-(t_{n-1}-t_{k-1})^{2-\alpha} -(t_n-t_k)^{2-\alpha}+(t_{n-1}-t_k)^{2-\alpha}}{\tau_k}
  M_h (U^k_{h\tau} - U^{k-1}_{h\tau}),
\end{align*}
the mass matrix
\begin{eqnarray}\label{chp-02-006}
  M_h = \frac{h}{6}\left(
   \begin{array}{ccccc}
    4 & 1 &   &  &  \\
    1 & 4 & 1 &  &  \\
      & \ddots & \ddots & \ddots &  \\
      &  & 1 & 4 & 1 \\
      &  &   & 1 & 4
   \end{array}
  \right)_{(M-1)\times(M-1)},
\end{eqnarray}
the stiffness matrix $A_h^\beta = (a^h_{i,j})_{(M-1)\times(M-1)}$ with its entries
\begin{align}\label{chp-02-001}
  \left \{
    \begin{aligned}
      & a^h_{i,i} = \frac{h^{1-2\beta}(2^{4-2\beta}-8)}{2\cos(\beta\pi)\Gamma(4-2\beta)},~~~\quad\qquad\qquad\qquad\qquad\qquad\qquad~~~~ i=1,\cdots,M-1\\
      & a^h_{j,j+1}=a^h_{j+1,j}=\frac{h^{1-2\beta}(3^{3-2\beta}-2^{5-2\beta}+7)}{2\cos(\beta\pi)\Gamma(4-2\beta)},~\quad\qquad\qquad\quad\quad j=1,\cdots,M-2\\
      & a^h_{k,k+l}= a^h_{k+l,k} = \frac{h^{1-2\beta}}{2\cos(\beta\pi)\Gamma(4-2\beta)}[(l+2)^{3-2\beta}\\
      & \quad\quad\quad-4(l+1)^{3-2\beta}+6l^{3-2\beta}-4(l-1)^{3-2\beta}+(l-2)^{3-2\beta}],~k=1,\cdots,M-l-1
    \end{aligned}
  \right.,
\end{align}
the vector 
\begin{eqnarray*}
  F^n_{h\tau}=(f^n_1,f^n_2,\cdots,f^n_{M-1})^T,~f^n_l=(f,\phi_l\times 1)_{Q_n},~l=1,\cdots,M-1
\end{eqnarray*}
and the fully FE approximations
\begin{eqnarray*}
  U^k_{h\tau}=(u^k_1,u^k_2,\cdots,u^k_{M-1})^T,~u^0_j=\psi_{0,I}(x_j),~
  u^k_j=u_h^k(x_j),~k=1,\cdots,n,~j=1,\cdots,M-1.
\end{eqnarray*}

\begin{Remark}
  \eqref{chp-02-009} is reduced via dividing both sides of \eqref{chp-02-004} by the factor
  $\tau^{1-\alpha}_n / \Gamma(3-\alpha)$, in case of the severe loss in convergence of the fully discrete FE scheme.
\end{Remark}

Next, a number of 
characterizations are established regarding $A_h^\beta$ just defined by \eqref{chp-02-001}.

\begin{Theorem}\label{thm-01}
  The stiffness matrix $A_h^\beta$ is symmetric and satisfies
\begin{enumerate}
  \item $a^h_{i,i} > 0$ for $i=1,\cdots,M-1$;
  \item $a^h_{i,j} < 0$ for $i\ne j$, $i,j=1,\cdots,M-1$;
  \item $\sum_{j=1}^{M-1} a^h_{i,j}>0$ for $i=1,\cdots,M-1$;
  \item The following relation holds for the particular case when $h\le1/7$
\begin{align*}
  \sum\limits_{j=1}^{M-1} a^h_{i,j}\ge \left \{
    \begin{aligned}
    & -\frac{h^{1-2\beta}(4-2^{3-2\beta})}{2\cos(\beta\pi)\Gamma(4-2\beta)},~i=1,M-1 \\
    & -\frac{2^{2\beta}h(2\beta-1)}{\cos(\beta\pi)\Gamma(2-2\beta)},~~~i=2,\cdots,M-2
    \end{aligned}
  \right.;
\end{align*}
  \item $A_h^\beta$ is an M-matrix.
\end{enumerate}
\end{Theorem}

\begin{Proof}
  The symmetric property of $A_h^\beta$ is an obvious fact by \eqref{chp-02-001}.
  Since $\beta\in(1/2,1)$, then $4-2\beta<3$ and $\cos(\beta\pi)<0$, which give immediately $a^h_{ii} > 0$, $i=1,\cdots,M-1$.
  This proves the first part of the theorem.
  The second part is an immediate consequence of the facts that on the interval $\beta\in(1/2,1)$,
  $f(\beta)=3^{3-2\beta}-2^{5-2\beta}+7$ is a strictly increasing function, and the bivariate function
\begin{eqnarray*}
  f_\beta(l)=(l+2)^{3-2\beta}-4(l+1)^{3-2\beta}+6l^{3-2\beta}-4(l-1)^{3-2\beta}+(l-2)^{3-2\beta} > 0,~2\le l\le M-2.
\end{eqnarray*}
  In fact, it is evident that
\begin{eqnarray*}
  1.5^{2\beta} > 1.5 > \frac{3^3}{2^5}\frac{\ln3}{\ln2} \Rightarrow f'(\beta)=-2\ln3 \cdot 3^{3-2\beta} + 2\ln2 \cdot 2^{5-2\beta} > 0,
\end{eqnarray*}
  and
\begin{eqnarray*}
  f_\beta(l)=h^{2\beta-3}[g(x_{l+2})-4g(x_{l+1})+6g(x_l)-4g(x_{l-1})+g(x_{l-2})] > 0
\end{eqnarray*}
  using Taylor's expansion with
\begin{eqnarray*}
  (\frac{l}{l+1})^{2+2\beta}-(\frac{l}{l-1})^{2+2\beta} > -\frac{30l}{2\beta+1} \Rightarrow
  l^{-1-2\beta}+\frac{2\beta+1}{30}[{(l+1)}^{-2-2\beta}-{(l-1)}^{-2-2\beta}] > 0,
\end{eqnarray*}
  where $g(x)=(x-a)^{3-2\beta}$ and $x_l=lh+a$.

  To prove the third part, use
\begin{eqnarray*}
  ({}_xD^{\beta}_L\phi_i,{}_xD^{\beta}_R\phi_j)_{L^2(\Omega)} = -({}_xD_L^{2\beta-1}\phi_i,\frac{d\phi_j}{dx})_{L^2(\Omega)}
\end{eqnarray*}
  and
\begin{eqnarray*}
  \tilde{\phi}:=\sum\limits_{j=1}^{M-1}\phi_j=1-\phi_0-\phi_M
\end{eqnarray*}
  to obtain the relation
\begin{eqnarray*}
  \sum\limits_{j=1}^{M-1} a^h_{i,j}=-\frac{({}_xD_L^{2\beta-1}\tilde{\phi},1)_{\Omega_i}-({}_xD_L^{2\beta-1}\tilde{\phi},1)_{\Omega_{i+1}}
  +({}_xD_L^{2\beta-1}\phi_i,1)_{\Omega_1}-({}_xD_L^{2\beta-1}\phi_i,1)_{\Omega_{M_h}}}{2h\cos(\beta\pi)},
\end{eqnarray*}
  where
\begin{align*}
  &{}_xD^{2\beta-1}_L\tilde{\phi}=\left \{
    \begin{aligned}
    & \frac{(x-a)^{2-2\beta}}{h\Gamma(3-2\beta)},~\qquad\qquad\qquad\qquad\qquad\qquad\qquad\quad ~a < x < x_1 \\
    & \frac{(x-a)^{2-2\beta}-(x-x_1)^{2-2\beta}}{h\Gamma(3-2\beta)},~\qquad\qquad\qquad\qquad ~~x_1 < x < x_{M_h-1} \\
    & \frac{(x-a)^{2-2\beta}-(x-x_1)^{2-2\beta}-(x-x_{M_h-1})^{2-2\beta}}{h\Gamma(3-2\beta)},~x_{M_h-1} < x < x_{M_h}
    \end{aligned}
  \right.
\end{align*}
  and
\begin{align*}
  {}_xD_L^{2\beta-1}\phi_i(x)=\left \{
    \begin{aligned}
    & 0,\qquad\qquad\qquad\qquad\qquad\qquad\qquad\qquad\qquad\qquad\quad ~x < x_{i-1} \\
    & \frac{(x-x_{i-1})^{2-2\beta}}{h\Gamma(3-2\beta)},\qquad\qquad\qquad\qquad\qquad\qquad\qquad\quad x_{i-1} < x < x_i \\
    & \frac{(x-x_{i-1})^{2-2\beta}-2(x-x_i)^{2-2\beta}}{h\Gamma(3-2\beta)},\qquad\qquad\qquad\qquad x_i < x < x_{i+1} \\
    & \frac{(x-x_{i-1})^{2-2\beta}-2(x-x_i)^{2-2\beta}+(x-x_{i+1})^{2-2\beta}}{h\Gamma(3-2\beta)},~x > x_{i+1}
    \end{aligned}
  \right..
\end{align*}
  Assume that $\Omega=(0,1)$ without loss of generality, one can easily derive
\begin{align*}
  &\sum\limits_{j=1}^{M-1} a^h_{i,j}=-\frac{(4-2^{3-2\beta})h^{3-2\beta}-1+3(1-h)^{3-2\beta}
  -3(1-2h)^{3-2\beta}+(1-3h)^{3-2\beta}}{2\cos(\beta\pi)h^2\Gamma(4-2\beta)},~i=1,M-1,\\
  & \sum\limits_{j=1}^{M-1} a^h_{i,j}=-\frac{3(ih)^{3-2\beta}-3[(i-1)h]^{3-2\beta}
  +[(i-2)h]^{3-2\beta}-[(i+1)h]^{3-2\beta}}{2\cos(\beta\pi)h^2\Gamma({4-2\beta})} \\
  & - \frac{3(1-ih)^{3-2\beta}-[1-(i-1)h]^{3-2\beta}-3[1-(i+1)h]^{3-2\beta}
  +[1-(i+2)h]^{3-2\beta}}{2\cos(\beta\pi)h^2\Gamma({4-2\beta})},~i=2,\cdots,M-2
\end{align*}
  and deduce $\sum_{j=1}^{M-1} a^h_{i,j}>0$ by Taylor's formula and $\beta\in(1/2,1)$.

  Another step to do in the proof is the result 4, which follows from
\begin{eqnarray*}
h\le1/7\Rightarrow7-4(1-\xi)^{-1-2\beta} > -\frac{1}{\beta h}
\Rightarrow-1+3(1-h)^{3-2\beta}-3(1-2h)^{3-2\beta}+(1-3h)^{3-2\beta}>0
\end{eqnarray*}
  for all $\xi\in(0,2h)$,
\begin{eqnarray*}
  (\frac{i-1}{i})^{3+2\beta} > \frac{\beta+1}{21(\beta+1)+30i}\Rightarrow
  \frac{3(ih)^{3-2\beta}-3(ih-h)^{3-2\beta}+(ih-2h)^{3-2\beta}-
  (ih+h)^{3-2\beta}}{h^3(3-2\beta)(2-2\beta)(2\beta-1)(ih)^{-2\beta}}>1
\end{eqnarray*}
  and
\begin{eqnarray*}
\frac{3(1-ih)^{3-2\beta}-[1-(i-1)h]^{3-2\beta}-3[1-(i+1)h]^{3-2\beta}
  +[1-(i+2)h]^{3-2\beta}}{h^3(3-2\beta)(2-2\beta)(2\beta-1)(1-ih)^{-2\beta}}>1
\end{eqnarray*}
  for $i=2,\cdots,M-2$, together with the inequality $(ih)^{-2\beta}+(1-ih)^{-2\beta}\ge 2^{1+2\beta}$.

  Finally, according to properties 1 and 2, the result 5 will be proved by showing that
  $(A_h^\beta)^{-1}$ is nonnegative, which can be easily proved by contradiction with property 3.
\end{Proof}

Observe from \eqref{chp-02-006}-\eqref{chp-02-001} that $M_h$ and $A_h^\beta$
are both symmetric Toeplitz matrices independent of any time terms.
The under-mentioned corollaries are natural consequences of Theorem \ref{thm-01}.

\begin{Corollary}\label{cor-01}
  The coefficient matrix $C^n_{h\tau}$ is a symmetric Toeplitz matrix. Furthermore,
  it will be independent of time level $n$ if the temporal discretization is also uniform.
\end{Corollary}

\begin{Corollary}\label{cor-02}
  The coefficient matrix $C^n_{h\tau}$ is an M-matrix, if and only if
\begin{eqnarray}\label{chp-03-015}
  \frac{\tau_n^\alpha}{h^{2\beta}} > -\frac{2\cos(\beta\pi)\Gamma(4-2\beta)}{3\Gamma(3-\alpha)(3^{3-2\beta}-2^{5-2\beta}+7)}.
\end{eqnarray}
\end{Corollary}

\begin{Proof}
  This result will follow from Theorem \ref{thm-01}, if we can show that
\begin{eqnarray*}
  \frac{h}{6}+\frac{\Gamma(3-\alpha)}{2}\tau^\alpha_n
  \frac{h^{1-2\beta}(3^{3-2\beta}-2^{5-2\beta}+7)}{2\cos(\beta\pi)\Gamma(4-2\beta)}<0,
\end{eqnarray*}
  which is an immediate application of the condition \eqref{chp-03-015}.
\end{Proof}

\subsection{Numerical experiments and the saturation error order}

\begin{Example}\label{exm-01}
  Consider \eqref{chp-01-01}-\eqref{chp-01-03} with $\Omega=(0,1)$, $T=1$, $\psi_0(x)=0$ and
\begin{align*}
  & f(x,t) =\frac{\Gamma(3-\alpha)}{\Gamma(3-2\alpha)}t^{2-2\alpha}x^2(1-x)^2 +
  \frac{t^{2-\alpha}}{\cos(\beta\pi)}\Big{[}\frac{x^{2-2\beta}+(1-x)^{2-2\beta}}{\Gamma(3-2\beta)}-\\
  &~~~~~~~~~~~~~~~~~~~~~~~\frac{6x^{3-2\beta}+6(1-x)^{3-2\beta}}{\Gamma(4-2\beta)} +
  \frac{12x^{4-2\beta}+12(1-x)^{4-2\beta}}{\Gamma(5-2\beta)}\Big{]}.
\end{align*}
\end{Example}

The exact solution is $u(x,t)=t^{2-\alpha}x^2(1-x)^2$. In the case of uniform temporal and spatial meshes,
Tables \ref{ctp-02-01} and \ref{ctp-02-02}
present errors $\|e\|_0 := \|u(\cdot,1)-u_{h\tau}(\cdot,1)\|_{L^2(\Omega)}$ and convergence rates.


\begin{table}[htbp]
\footnotesize
\centering\caption{Error results and convergence rates in spatial direction with $h=\tau$.}\label{ctp-02-01}\vskip 0.1cm
\begin{tabular}{|c|c|c|c|c|c|c|c|c|c|c|c|c|c|c|c|c|c|c|}\hline
\multirow{3}{*}{$N$}&\multicolumn{6}{|c|}{$\beta=0.6$}&\multicolumn{6}{|c|}{$\beta=0.8$} \\ \cline{2-13}
~&\multicolumn{2}{|c|}{$\alpha=0.01$}&\multicolumn{2}{|c|}{$\alpha=0.50$}&\multicolumn{2}{|c|}{$\alpha=0.99$}&
\multicolumn{2}{|c|}{$\alpha=0.01$}&\multicolumn{2}{|c|}{$\alpha=0.50$}&\multicolumn{2}{|c|}{$\alpha=0.99$} \\ \cline{2-13}
~& $\|e\|_0$ & rate & $\|e\|_0$ & rate & $\|e\|_0$ & rate & $\|e\|_0$ & rate & $\|e\|_0$ & rate & $\|e\|_0$ & rate \\\hline

  8 & 8.94E-4 &   -  & 8.78E-4 &   -  & 8.88E-4 &   -  & 9.73E-4 &   -  & 9.56E-4 &  -   & 9.72E-4 &   -  \\ \hline
 16 & 2.02E-4 & 2.15 & 1.98E-4 & 2.15 & 2.00E-4 & 2.15 & 2.34E-4 & 2.06 & 2.29E-4 & 2.06 & 2.33E-4 & 2.06 \\ \hline
 32 & 4.49E-5 & 2.17 & 4.42E-5 & 2.16 & 4.47E-5 & 2.16 & 5.51E-5 & 2.08 & 5.40E-5 & 2.08 & 5.49E-5 & 2.08 \\ \hline
 64 & 1.01E-5 & 2.15 & 1.00E-5 & 2.14 & 1.01E-5 & 2.14 & 1.30E-5 & 2.09 & 1.27E-5 & 2.09 & 1.29E-5 & 2.09 \\ \hline

\multirow{2}{*}{$N$}&\multicolumn{2}{|c|}{$\alpha=0.10$}
&\multicolumn{2}{|c|}{$\alpha=0.25$}&\multicolumn{2}{|c|}{$\alpha=0.75$} &\multicolumn{2}{|c|}{$\alpha=0.10$}
&\multicolumn{2}{|c|}{$\alpha=0.25$}&\multicolumn{2}{|c|}{$\alpha=0.75$} \\ \cline{2-13}
~& $\|e\|_0$ & rate & $\|e\|_0$ & rate & $\|e\|_0$ & rate & $\|e\|_0$ & rate & $\|e\|_0$ & rate & $\|e\|_0$ & rate \\\hline

  8 & 8.90E-4 &   -  & 8.80E-4 &  -   & 8.83E-4 &   -  & 9.72E-4 &   -  & 9.63E-4 &  -   & 9.64E-4 &   -  \\ \hline
 16 & 1.99E-4 & 2.16 & 1.97E-4 & 2.16 & 1.99E-4 & 2.15 & 2.31E-4 & 2.07 & 2.28E-4 & 2.08 & 2.31E-4 & 2.06 \\ \hline
 32 & 4.39E-5 & 2.18 & 4.40E-5 & 2.16 & 4.44E-5 & 2.16 & 5.39E-5 & 2.10 & 5.34E-5 & 2.09 & 5.44E-5 & 2.09 \\ \hline
 64 & 9.95E-6 & 2.14 & 1.00E-5 & 2.14 & 1.01E-5 & 2.14 & 1.25E-5 & 2.11 & 1.26E-5 & 2.08 & 1.28E-5 & 2.09 \\ \hline

\end{tabular}
\end{table}

\begin{table}[htbp]
\footnotesize
\centering\caption{Error results and convergence rates in spatial direction with $h=\sqrt{\tau}$.}\label{ctp-02-02}\vskip 0.1cm
\begin{tabular}{|c|c|c|c|c|c|c|c|c|c|c|c|c|c|c|c|c|c|c|}\hline
\multirow{3}{*}{$N$}&\multicolumn{6}{|c|}{$\beta=0.6$}&\multicolumn{6}{|c|}{$\beta=0.8$} \\ \cline{2-13}
~&\multicolumn{2}{|c|}{$\alpha=0.01$}&\multicolumn{2}{|c|}{$\alpha=0.50$}&\multicolumn{2}{|c|}{$\alpha=0.99$}&
\multicolumn{2}{|c|}{$\alpha=0.01$}&\multicolumn{2}{|c|}{$\alpha=0.50$}&\multicolumn{2}{|c|}{$\alpha=0.99$} \\ \cline{2-13}
~& $\|e\|_0$ & rate & $\|e\|_0$ & rate & $\|e\|_0$ & rate & $\|e\|_0$ & rate & $\|e\|_0$ & rate & $\|e\|_0$ & rate \\\hline

 16  & 3.64E-3 &   -  & 3.63E-3 &   -  & 3.63E-3 &   -  & 3.76E-3 &   -  & 3.75E-3 &   -  & 3.75E-3 &  -   \\ \hline
 64  & 8.94E-4 & 1.01 & 8.87E-4 & 1.02 & 8.88E-4 & 1.02 & 9.73E-4 & 0.97 & 9.69E-4 & 0.98 & 9.72E-4 & 0.98 \\ \hline
 256 & 2.02E-4 & 1.07 & 2.00E-4 & 1.08 & 2.00E-4 & 1.08 & 2.34E-4 & 1.03 & 2.32E-4 & 1.03 & 2.33E-4 & 1.03 \\ \hline

\multirow{2}{*}{$N$}&\multicolumn{2}{|c|}{$\alpha=0.10$}
&\multicolumn{2}{|c|}{$\alpha=0.25$}&\multicolumn{2}{|c|}{$\alpha=0.75$} &\multicolumn{2}{|c|}{$\alpha=0.10$}
&\multicolumn{2}{|c|}{$\alpha=0.25$}&\multicolumn{2}{|c|}{$\alpha=0.75$} \\ \cline{2-13}
~& $\|e\|_0$ & rate & $\|e\|_0$ & rate & $\|e\|_0$ & rate & $\|e\|_0$ & rate & $\|e\|_0$ & rate & $\|e\|_0$ & rate \\\hline

 16  & 3.64E-3 &   -  & 3.63E-3 &   -  & 3.63E-3 &   -  & 3.75E-3 &   -  & 3.75E-3 &   -  & 3.75E-3 &   -  \\ \hline
 64  & 8.90E-4 & 1.02 & 8.87E-4 & 1.02 & 8.86E-4 & 1.02 & 9.72E-4 & 0.98 & 9.71E-4 & 0.98 & 9.69E-4 & 0.98 \\ \hline
 256 & 2.01E-4 & 1.07 & 2.01E-4 & 1.07 & 2.00E-4 & 1.08 & 2.33E-4 & 1.03 & 2.33E-4 & 1.03 & 2.32E-4 & 1.03 \\ \hline

\end{tabular}
\end{table}

From Tables \ref{ctp-02-01} and \ref{ctp-02-02}, we can obtain that
the fully FE solution $u_{h\tau}$ achieves the saturation error order $\mathcal{O}(\tau^2+h^2)$ under $\|\cdot\|_0$ norm.

\begin{figure}[htp]
  \centerline{
  \includegraphics[scale=0.4]{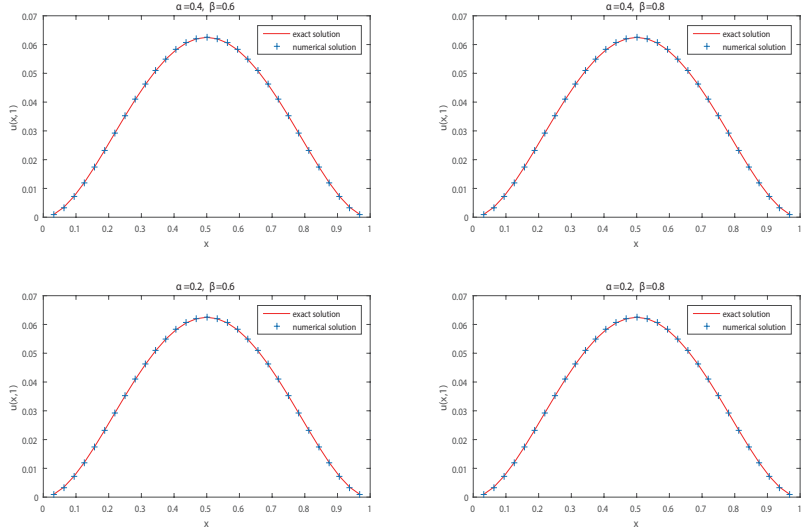}}
  \caption{Illustration for comparisons on exact solutions and numerical solutions with $t=1$ and $h=\tau=1/32$.}
  \label{fig-01}
\end{figure}

Fig. \ref{fig-01} illustrates the comparisons of exact solutions and numerical solutions
of $\alpha=0.2$, $0.4$ and $\beta=0.6$, $0.8$ with $t=1$ and $h=\tau=1/32$.

\section{Condition number estimation}

This section is devoted to deriving the condition number estimation on the coefficient matrix of \eqref{chp-02-009}
in uniform temporal and and spatial discretizations.

\begin{Theorem}\label{clp-02-04}
  For the linear system \eqref{chp-02-009}, we have
\begin{eqnarray}\label{chp-02-007}
  \kappa(C^n_{h\tau})= 1+\mathcal{O}(\tau^\alpha h^{-2\beta}).
\end{eqnarray}
\end{Theorem}

\begin{Proof}
  Let $C_\alpha = \Gamma(3-\alpha)/2$, we divide our proof in three steps.
  First, it is trivially true that $C^n_{h\tau}$
  is spectrally equivalent to the matrix $I+C_\alpha\tau^\alpha M_h^{-\frac{1}{2}}A_h^\beta M_h^{-\frac{1}{2}}$, i.e.
\begin{eqnarray}\label{chp-03-012}
 \kappa(C^n_{h\tau})\simeq\kappa\Big{(}I+C_\alpha\tau^\alpha M_h^{-\frac{1}{2}}A_h^\beta M_h^{-\frac{1}{2}}\Big{)}.
\end{eqnarray}

  The next thing to do in the proof is to verify
\begin{eqnarray}\label{chp-03-010}
 \lambda_{\min}(M_h^{-\frac{1}{2}}A_h^\beta M_h^{-\frac{1}{2}})\gtrsim 1,~
 \lambda_{\max}(M_h^{-\frac{1}{2}}A_h^\beta M_h^{-\frac{1}{2}})\lesssim h^{-2\beta},
\end{eqnarray}
  which is equivalent to
\begin{eqnarray}\label{chp-03-009}
  (\vec{v}_h,\vec{v}_h)\lesssim (M_h^{-\frac{1}{2}}A_h^\beta M_h^{-\frac{1}{2}}\vec{v}_h,\vec{v}_h)
  \lesssim h^{-2\beta}(\vec{v}_h,\vec{v}_h),~\forall \vec{v}_h\in \mathbb{R}^{M-1}.
\end{eqnarray}
  Set $\vec{u}_h = M_h^{-\frac{1}{2}}\vec{v}_h := (u^h_1,\cdots,u^h_{M-1})^T$, rewrite \eqref{chp-03-009} as
  $(M_h\vec{u}_h,\vec{u}_h)\lesssim (A_h^\beta \vec{u}_h,\vec{u}_h)\lesssim h^{-2\beta}(M_h\vec{u}_h,\vec{u}_h)$.
  It is sufficient to verify that
  $(M_h\vec{u}_h,\vec{u}_h)\simeq h(\vec{u}_h,\vec{u}_h)$.
  It follows by \eqref{chp-02-006} and the Cauchy-Schwarz inequality that
\begin{eqnarray*}
  \frac{h}{3}(\vec{u}_h, \vec{u}_h)\le(M_h \vec{u}_h, \vec{u}_h)=h\Big{[}\frac{2}{3}\sum_{l=1}^{M-1}(u^h_l)^2
  +\frac{1}{3}\sum_{l=1}^{M-2}u^h_lu^h_{l+1}\Big{]}\le h(\vec{u}_h, \vec{u}_h).
\end{eqnarray*}
  Thus \eqref{chp-03-010} will follow if we can show that
  $h(\vec{u}_h,\vec{u}_h)\lesssim (A_h^\beta \vec{u}_h,\vec{u}_h)\lesssim h^{1-2\beta}(\vec{u}_h,\vec{u}_h).$
  We start by showing the second inequality. Utilizing Theorem \ref{thm-01} and the Cauchy-Schwarz inequality, we arrive at
\begin{eqnarray*}
 (A_h^\beta \vec{u}_h, \vec{u}_h) &\le& \sum_{i=1}^{M-1} a^h_{i,i} (u^h_i)^2 -
 \frac{1}{2} \sum_{i=1}^{M-1} \sum_{j\ne i} a^h_{i,j} [(u^h_i)^2 + (u^h_j)^2] \\
 &=& \sum_{i=1}^{M-1} a^h_{i,i} (u^h_i)^2 - \sum_{i=1}^{M-1} \sum_{j=i+1}^{M-1} a^h_{i,j} (u^h_i)^2
  - \sum_{i=1}^{M-1} \sum_{j=i+1}^{M-1} a^h_{i,j} (u^h_j)^2 \\
 &=& \sum_{i=1}^{M-1} (u^h_i)^2 \Big{[}a^h_{i,i} - \sum_{j=i+1}^{M-1} a^h_{i,j} - \sum_{j=1}^{i-1} a^h_{i,j}\Big{]} \\
 & \le & 2 a^h_{1,1} (\vec{u}_h, \vec{u}_h) = \frac{(2^{4-2\beta}-8)}{\cos(\beta\pi)\Gamma(4-2\beta)} h^{1-2\beta} (\vec{u}_h, \vec{u}_h),
\end{eqnarray*}
  which proves the second inequality. To prove the left inequality, set $u_h:= \Phi_h\vec{u}_h$, we rewrite it as
\begin{eqnarray*}
  \frac{1}{\cos(\beta\pi)}({}_xD^{\beta}_Lu_h,{}_xD^{\beta}_Ru_h)_{L^2(\Omega)}=(A_h^\beta \vec{u}_h,\vec{u}_h)
  \gtrsim h(\vec{u}_h,\vec{u}_h)\simeq (M_h\vec{u}_h,\vec{u}_h)= (u_h, u_h)_{L^2(\Omega)}
\end{eqnarray*}
  which can be deduced by \eqref{chp-03-002}-\eqref{chp-03-003}, where $\Phi_h = (\phi_1,\cdots,\phi_{M-1})$.

  Finally, we have to show that
\begin{eqnarray*}
  \kappa\Big{(}I+C_\alpha\tau^\alpha M_h^{-\frac{1}{2}}A_h^\beta M_h^{-\frac{1}{2}}\Big{)}&=&
  \frac{\lambda_{\max}\Big{(}I+C_\alpha\tau^\alpha M_h^{-\frac{1}{2}}A_h^\beta
  M_h^{-\frac{1}{2}}\Big{)}}{\lambda_{\min}\Big{(}I+C_\alpha\tau^\alpha M_h^{-\frac{1}{2}}A_h^\beta M_h^{-\frac{1}{2}}\Big{)}}\\
  &=&\frac{1+C_\alpha\tau^\alpha \lambda_{\max}(M_h^{-\frac{1}{2}}A_h^\beta
  M_h^{-\frac{1}{2}})}{1+C_\alpha\tau^\alpha \lambda_{\min}(M_h^{-\frac{1}{2}}A_h^\beta M_h^{-\frac{1}{2}})} \\
  &\le& 1+C_\alpha\frac{(2^{4-2\beta}-8)}{\cos(\beta\pi)\Gamma(4-2\beta)}\tau^\alpha h^{-2\beta}.
\end{eqnarray*}
  This completes the proof based on the spectral equivalence relation \eqref{chp-03-012}.
\end{Proof}

\begin{Remark}
  The estimation \eqref{chp-02-007} is compatible with the correlative result $1+\mathcal{O}(\tau h^{-2})$ of
  integer order parabolic differential equations.
\end{Remark}

An important particular case of Theorem \ref{clp-02-04} is singled out in the following corollary.

\begin{Corollary}\label{cor-03}
  Let $\tau$ be proportional to $h^\mu$ with $\mu\alpha\ge2\beta$. Then
\begin{eqnarray}\label{chp-03-013}
  \kappa(C^n_{h\tau})=\mathcal{O}(1).
\end{eqnarray}
\end{Corollary}

In what follows, we examine the correctness of \eqref{chp-02-007} concerning Example \ref{exm-01} with typical
$\alpha$ and $\beta$ for three specific cases: $\tau = h$, $\tau = h^2$ and $\tau$ is fixed (doesn't change along with $h$).
In under-mentioned tables, $\lambda_{\min}$ and $\lambda_{\max}$ respectively indicate the smallest and largest eigenvalues,
$\kappa$ represents the condition number and \emph{ratio} is the quotient of the condition number in fine grid divided by that in coarse grid.

\begin{table}[htbp]
\footnotesize
\centering\caption{The smallest and largest eigenvalues and condition numbers with $\tau=h$.}\label{ctp-03-01}\vskip 0.1cm
\begin{tabular}{|c|c|c|c|c|c|c|c|c|c|}\hline
\multirow{2}{*}{$\alpha$}&\multirow{2}{*}{$M$}&\multicolumn{4}{|c|}{$\beta=0.6$}
&\multicolumn{4}{|c|}{$\beta=0.8$} \\ \cline{3-10} ~&~&
$\lambda_{\min}$ & $\lambda_{\max}$ & $\kappa$ & \emph{ratio} &
$\lambda_{\min}$ & $\lambda_{\max}$ & $\kappa$ & \emph{ratio} \\\hline

\multirow{4}{*}{0.99} & 8 & 1.45E-1 & 1.97E-1 & 1.36E+0 &  -   & 1.63E-1 & 5.47E-1 & 3.35E+0 &  -   \\ \cline{2-10}
& 16 & 6.81E-2 & 1.09E-1 & 1.60E+0 & 1.18 & 7.26E-2 & 4.11E-1 & 5.66E+0 & 1.69 \\ \cline{2-10}
& 32 & 3.27E-2 & 6.03E-2 & 1.84E+0 & 1.15 & 3.39E-2 & 3.09E-1 & 9.12E+0 & 1.61 \\ \cline{2-10}
& 64 & 1.60E-2 & 3.34E-2 & 2.09E+0 & 1.13 & 1.63E-2 & 2.33E-1 & 1.43E+1 & 1.57 \\ \hline

\multirow{4}{*}{0.5} & 8 & 2.09E-1 & 5.41E-1 & 2.59E+0 &  -   & 2.74E-1 & 1.89E+0 & 6.89E+0 &  -   \\ \cline{2-10}
& 16 & 9.31E-2 & 4.25E-1 & 4.56E+0 & 1.76 & 1.16E-1 & 2.03E+0 & 1.74E+1 & 2.53 \\ \cline{2-10}
& 32 & 4.22E-2 & 3.37E-1 & 8.00E+0 & 1.75 & 5.05E-2 & 2.17E+0 & 4.30E+1 & 2.47 \\ \cline{2-10}
& 64 & 1.95E-2 & 2.70E-1 & 1.39E+1 & 1.73 & 2.24E-2 & 2.32E+0 & 1.04E+2 & 2.41 \\ \hline

\multirow{4}{*}{0.01} & 8 & 4.81E-1 & 2.08E+0 & 4.32E+0 &  -   & 7.50E-1 & 7.65E+0 & 1.02E+1 &  -   \\ \cline{2-10}
& 16 & 2.42E-1 & 2.35E+0 & 9.69E+0 & 2.24 & 3.77E-1 & 1.17E+1 & 3.09E+1 & 3.03 \\ \cline{2-10}
& 32 & 1.21E-1 & 2.66E+0 & 2.21E+1 & 2.28 & 1.88E-1 & 1.76E+1 & 9.36E+1 & 3.03 \\ \cline{2-10}
& 64 & 6.00E-2 & 3.03E+0 & 5.06E+1 & 2.29 & 9.35E-2 & 2.65E+1 & 2.83E+2 & 3.03 \\ \hline

\end{tabular}
\end{table}

\begin{table}[htbp]
\footnotesize
\centering\caption{The smallest and largest eigenvalues and condition numbers with $\tau=h^2$.}\label{ctp-03-02}\vskip 0.1cm
\begin{tabular}{|c|c|c|c|c|c|c|c|c|c|}\hline
\multirow{2}{*}{$\alpha$}&\multirow{2}{*}{$M$}&\multicolumn{4}{|c|}{$\beta=0.6$}
&\multicolumn{4}{|c|}{$\beta=0.8$} \\ \cline{3-10} ~&~&
$\lambda_{\min}$ & $\lambda_{\max}$ & $\kappa$ & \emph{ratio} &
$\lambda_{\min}$ & $\lambda_{\max}$ & $\kappa$ & \emph{ratio} \\\hline

\multirow{4}{*}{0.5} & 8 & 1.52E-1 & 2.33E-1 & 1.53E+0 &  -   &  1.76E-1 & 6.96E-1 & 3.97E+0 &  -   \\ \cline{2-10}
& 16 & 6.98E-2 & 1.30E-1 & 1.86E+0 & 1.21 &  7.57E-2 & 5.23E-1 & 6.91E+0 & 1.74 \\ \cline{2-10}
& 32 & 3.31E-2 & 7.20E-2 & 2.17E+0 & 1.17 &  3.46E-2 & 3.92E-1 & 1.13E+1 & 1.64 \\ \cline{2-10}
& 64 & 1.61E-2 & 4.01E-2 & 2.49E+0 & 1.15 &  1.65E-2 & 2.95E-1 & 1.79E+1 & 1.58 \\ \hline

\multirow{4}{*}{0.01} & 8 & 4.74E-1 & 2.04E+0 & 4.30E+0 &  -   & 7.37E-1 & 7.49E+0 & 1.02E+1 &  -   \\ \cline{2-10}
& 16 & 2.37E-1 & 2.28E+0 & 9.62E+0 & 2.24 & 3.69E-1 & 1.13E+1 & 3.08E+1 & 3.02 \\ \cline{2-10}
& 32 & 1.18E-1 & 2.57E+0 & 2.19E+1 & 2.27 & 1.83E-1 & 1.70E+1 & 9.30E+1 & 3.03 \\ \cline{2-10}
& 64 & 5.82E-2 & 2.91E+0 & 5.00E+1 & 2.29 & 9.03E-2 & 2.54E+1 & 2.81E+2 & 3.03 \\ \hline

\multirow{4}{*}{$\beta$} & 8 & 1.40E-1 & 1.75E-1 & 1.24E+0 &  -   & 1.35E-1 & 2.00E-1 & 1.49E+0 &  -   \\ \cline{2-10}
& 16 & 6.62E-2 & 8.84E-2 & 1.34E+0 & 1.07 & 6.42E-2 & 1.00E-1 & 1.56E+0 & 1.05 \\ \cline{2-10}
& 32 & 3.21E-2 & 4.44E-2 & 1.38E+0 & 1.04 & 3.16E-2 & 5.00E-2 & 1.59E+0 & 1.02 \\ \cline{2-10}
& 64 & 1.58E-2 & 2.22E-2 & 1.40E+0 & 1.02 & 1.57E-2 & 2.50E-2 & 1.60E+0 & 1.01 \\ \hline

\end{tabular}
\end{table}

\begin{table}[htbp]
\footnotesize
\centering\caption{The smallest and largest eigenvalues and the condition number with $\tau=1/32$.}\label{ctp-03-03}\vskip 0.1cm
\begin{tabular}{|c|c|c|c|c|c|c|c|c|c|}\hline
\multirow{2}{*}{$\alpha$}&\multirow{2}{*}{$M$}&\multicolumn{4}{|c|}{$\beta=0.6$}
&\multicolumn{4}{|c|}{$\beta=0.8$} \\ \cline{3-10} ~&~&
$\lambda_{\min}$ & $\lambda_{\max}$ & $\kappa$ & \emph{ratio} &
$\lambda_{\min}$ & $\lambda_{\max}$ & $\kappa$ & \emph{ratio} \\\hline

\multirow{4}{*}{0.99} & 64 & 1.64E-2 & 5.83E-2 & 3.56E+0 & 1.93 & 1.69E-2 & 4.58E-1 & 2.70E+1 & 2.96 \\ \cline{2-10}
& 128 & 8.19E-3 & 6.24E-2 & 7.62E+0 & 2.14 & 8.48E-3 & 6.89E-1 & 8.13E+1 & 3.01 \\ \cline{2-10}
& 256 & 4.10E-3 & 6.97E-2 & 1.70E+1 & 2.23 & 4.24E-3 & 1.04E+0 & 2.46E+2 & 3.02 \\ \cline{2-10}
& 512 & 2.05E-3 & 7.91E-2 & 3.86E+1 & 2.27 & 2.12E-3 & 1.58E+0 & 7.45E+2 & 3.03 \\ \hline

\multirow{4}{*}{0.5} & 64 & 2.11E-2 & 3.80E-1 & 1.80E+1 & 2.25 & 2.52E-2 & 3.28E+0 & 1.30E+2 & 3.02 \\ \cline{2-10}
& 128 & 1.06E-2 & 4.33E-1 & 4.10E+1 & 2.28 & 1.26E-2 & 4.97E+0 & 3.94E+2 & 3.03 \\ \cline{2-10}
& 256 & 5.27E-3 & 4.95E-1 & 9.39E+1 & 2.29 & 6.31E-3 & 7.53E+0 & 1.19E+3 & 3.03 \\ \cline{2-10}
& 512 & 2.64E-3 & 5.68E-1 & 2.15E+2 & 2.29 & 3.16E-3 & 1.14E+1 & 3.62E+3 & 3.03 \\ \hline

\multirow{4}{*}{0.01} & 64 & 6.03E-2 & 3.05E+0 & 5.07E+1 & 2.29 & 9.40E-2 & 2.67E+1 & 2.84E+2 & 3.03 \\ \cline{2-10}
& 128 & 3.01E-2 & 3.50E+0 & 1.16E+2 & 2.30 & 4.70E-2 & 4.05E+1 & 8.61E+2 & 3.03 \\ \cline{2-10}
& 256 & 1.51E-2 & 4.02E+0 & 2.67E+2 & 2.30 & 2.35E-2 & 6.13E+1 & 2.61E+3 & 3.03 \\ \cline{2-10}
& 512 & 7.53E-3 & 4.62E+0 & 6.14E+2 & 2.30 & 1.17E-2 & 9.29E+1 & 7.91E+3 & 3.03 \\ \hline \hline

\multirow{2}{*}{$\beta$}&\multirow{2}{*}{$M$}&\multicolumn{4}{|c|}{$\alpha=0.01$}
&\multicolumn{4}{|c|}{$\alpha=0.99$} \\ \cline{3-10} ~&~&
$\lambda_{\min}$ & $\lambda_{\max}$ & $\kappa$ & \emph{ratio} &
$\lambda_{\min}$ & $\lambda_{\max}$ & $\kappa$ & \emph{ratio} \\\hline

\multirow{4}{*}{0.999} & 64 & 1.63E-1 & 2.42E+2 & 1.49E+3 & 4.00 & 1.81E-2 & 4.11E+0 & 2.27E+2 & 3.98 \\ \cline{2-10}
& 128 & 8.14E-2 & 4.84E+2 & 5.95E+3 & 4.00 & 9.06E-3 & 8.21E+0 & 9.07E+2 & 3.99 \\ \cline{2-10}
& 256 & 4.07E-2 & 9.66E+2 & 2.38E+4 & 3.99 & 4.53E-3 & 1.64E+1 & 3.62E+3 & 3.99 \\ \cline{2-10}
& 512 & 2.03E-2 & 1.93E+3 & 9.49E+4 & 3.99 & 2.27E-3 & 3.28E+1 & 1.45E+4 & 3.99 \\ \hline

\end{tabular}
\end{table}

It is observed from Tables \ref{ctp-03-01}-\ref{ctp-03-03} that numerical results are in good agreement with our theoretical estimation.

\section{AMG's convergence analysis and an adaptive AMG method}

Within the section, involving FFT to perform Toeplitz matrix-vector multiplications,
we introduce the so-called Ruge-St{\"u}ben or classical AMG method \cite{r-001}
with low algorithmic complexity, fulfill its theoretical investigation,
and then propose an adaptive AMG method through Corollary \ref{cor-03}.

\begin{Algorithm}\label{alg-02} The classical AMG method for the linear system \eqref{chp-02-009}.
  \begin{description}

    \item[Step 1] Perform the Setup phase to the coefficient matrix $C^n_{h\tau}$.
    \begin{description}

      \item[1.1] Set the strength-of-connection tolerance $\theta$;

      \item[1.2] Build the ingredients required by a hierarchy of levels, coarsest to finest, including the grid transfer operator $P$.

    \end{description}

    \item[Step 2] Invoke the classical V($\varrho_1$,$\varrho_2$)-cycle to solve \eqref{chp-02-009} until convergence.
    Below is the description of two-grid V($\varrho_1$,$\varrho_2$)-cycle.
    \begin{description}

      \item[2.1] Do $\varrho_1$ pre-smoothing steps on \eqref{chp-02-009};

      \item[2.2] Compute and restrict the residual: $r^c=P^T(G^n_{h\tau} - C^n_{h\tau}U^n_{h\tau})$;

      \item[2.3] Solve the residual equation on coarse level: $(P^T C^n_{h\tau} P) e^c = r^c$;

      \item[2.4] Interpolation and correction: $U^n_{h\tau}=U^n_{h\tau}+Pe^c$;

      \item[2.5] Do $\varrho_2$ post-smoothing steps on \eqref{chp-02-009}.

    \end{description}

  \end{description}
\end{Algorithm}

\begin{Remark}
  In pre- and post-smoothing processes, damped-Jacobi iterative methods are favorable choices,
  which can maintain the low computational cost $\mathcal{O}(M\log M)$ calculated by FFT.
\end{Remark}

For theoretical investigations, we rewrite \eqref{chp-02-009} and
the grid transfer operator $P$ in block form regarding a given C/F splitting
\begin{eqnarray*}
C^n_{h\tau}U^n_{h\tau} =
\left(
\begin{array}{cc}
A_{FF} & A_{FC} \\
A_{CF} & A_{CC}
\end{array}
\right)
\left(
\begin{array}{c}
u_{F} \\
u_{C}
\end{array}
\right)
=\left(
\begin{array}{c}
f_{F} \\
f_{C}
\end{array}
\right) = G^n_{h\tau},~
P=\left(
\begin{array}{l}
I_{FC} \\
I_{CC}
\end{array}\right),
\end{eqnarray*}
and introduce the following inner products
\begin{eqnarray*}
(u_F,v_F)_{0,F} = (D_{FF}u_F,v_F),~(u,v)_1 =(C^n_{h\tau}u,v),~(u,v)_2 =(D^{-1}_{h\tau}C^n_{h\tau}u,C^n_{h\tau}v)
\end{eqnarray*}
with their associated norms $\|\cdot\|_{0,F}=\sqrt{(\cdot,\cdot)_{0,F}}$ and $\|\cdot\|_i=\sqrt{(\cdot,\cdot)_i}$ ($i=1,2$),
where $I_{CC}$ is the identity operator, $D_{FF}={\bf diag}(A_{FF})$ and $D_{h\tau}={\bf diag}(C^n_{h\tau})$.

For simplicity, we here denote $C^n_{h\tau}=(c_{ij})_{(M-1)\times(M-1)}$,
and only consider the two-grid V(0,1)-cycle, whose iteration matrix has the form
\begin{eqnarray*}
M_{h,H}=S[I-P(P^T C^n_{h\tau} P)^{-1}P^TC^n_{h\tau}],
\end{eqnarray*}
where $S$ is a relaxation operator usually chosen as damped-Jacobi or Gauss-Seidel iterative method.

Combining Corollary \ref{cor-02} and the two-level convergence theory in the work \cite{s-001}, leads to the following lemmas and theorem.

\begin{Lemma}\label{clp-03-01}
  Under the condition \eqref{chp-03-015}, for all $e_h\in \mathbb{R}^{M-1}$, damped-Jacobi and Gauss-Seidel relaxations satisfy the smoothing property
\begin{eqnarray}\label{chp-04-02}
  \|Se_h\|^2_{1}\leq \|e_h\|^2_1-\sigma_1\|e_h\|^2_2
\end{eqnarray}
  with $\sigma_1$ independent of $e_h$ and step sizes $h$ and $\tau_n$.
\end{Lemma}

\begin{Proof}
  On the strength of Theorem A.3.1 and A.3.2 in \cite{s-001},
  we produce that damped-Jacobi relaxation with parameter $0<\omega<2/\eta$ satisfies \eqref{chp-04-02} with
\begin{eqnarray*}
  \sigma_1=\omega(2-\omega\eta),
\end{eqnarray*}
  and Gauss-Seidel relaxation satisfies \eqref{chp-04-02} with
\begin{eqnarray*}
  \sigma_1=\frac{1}{(1+\gamma_-)(1+\gamma_+)},
\end{eqnarray*}
  both independent of $e_h$, where
\begin{eqnarray*}
  \eta \ge \rho(D^{-1}_{h\tau}C^n_{h\tau}),~
  \gamma_-=\max_i\Big{\{}\frac{1}{w_ic_{ii}}\sum_{j<i}w_j|c_{ij}|\Big{\}},~
  \gamma_+=\max_i\Big{\{}\frac{1}{w_ic_{ii}}\sum_{j>i}w_j|c_{ij}|\Big{\}},
\end{eqnarray*}
  and $w=(w_i)$ is an arbitrary positive vector with $C^n_{h\tau}w$ being also positive. 

  By exploiting \eqref{chp-02-008}-\eqref{chp-02-001}, the assumption \eqref{chp-03-015} and Corollary \ref{cor-02},
  we conclude that $C^n_{h\tau}$ is strictly diagonally dominant. Recall that $C^n_{h\tau}w$ is a positive vector,
  yield $\gamma_-<1$, $\gamma_+<1$ and
\begin{eqnarray}\label{chp-03-001}
  \rho(D^{-1}_{h\tau}C^n_{h\tau})\le|D^{-1}_{h\tau}C^n_{h\tau}|_w=\max_i\Big{\{}\frac{1}{w_i}\sum_jw_j\frac{|c_{ij}|}{c_{ii}}\Big{\}}<2,
\end{eqnarray}
  which implicitly mean that $\eta$, $\gamma_-$ and $\gamma_+$ can be chosen to be independent of $h$ and $\tau_n$, and complete the proof.
\end{Proof}

\begin{Remark}
  The inequality \eqref{chp-03-001} implies that there exists $\epsilon>0$ such that $\rho(D^{-1}_{h\tau}C^n_{h\tau})=2-3\epsilon$.
  Then $\eta=2-2\epsilon>\rho(D^{-1}_{h\tau}C^n_{h\tau})$ and hence the upper bound of parameter $\omega$: $2/\eta=1/(1-\epsilon)>1$,
  which suggests that Jacobi relaxation with $\omega=1$ is available in such a case.
\end{Remark}

\begin{Remark}
  For all symmetric M-matrices, $\sigma_1\le1/\eta<1$ holds for damped-Jacobi relaxation,
  while $\sigma_1\in(1/4,1)$ for Gauss-Seidel relaxation.
\end{Remark}

\begin{Lemma}\label{clp-03-02}
  Under the condition \eqref{chp-03-015} and a given C/F splitting,
  for all $e_h=(e_F^T,e_C^T)^T\in \mathbb{R}^{M-1}$,
  the direct interpolation $I_{FC}$ satisfies
\begin{eqnarray}\label{chp-04-03}
  \|e_F-I_{FC}e_C\|^2_{0,F} \le \sigma_2\|e_h\|^2_1
\end{eqnarray}
  with $\sigma_2$ independent of $e_h$, $h$ and $\tau_n$.
\end{Lemma}

\begin{Proof}
  According to Theorem A.4.3 in \cite{s-001}, $I_{FC}$ satisfies \eqref{chp-04-03} with
  $\sigma_2$ of the form regarding a given C/F splitting
\begin{eqnarray}\label{chp-04-05}
  \sigma_2 \ge \max_{i\in F} \Big{\{}\frac{\sum_{j\in N_i}c_{ij}}{\sum_{j\in C_i}c_{ij}}\Big{\}}
\end{eqnarray}
  independent of $e_h$, where $N_i=\{j\ne i:c_{ij}\ne0\}$, $C_i$ is the subset of $N_i$
  whose values will be used to interpolate at F-point $i$. As a result of \eqref{chp-04-05}
  and the fact that $c_{ij}$ ($j\in N_i$) are all negative, the following relation holds: $\sigma_2 > 1$.

  Notice here that the classical Ruge-St{\"u}ben based coarsening strategy generates
  at least one of points $i-1$ and $i+1$ to be C-points and strongly influence $i$
  --- viz. it retains $i-1\in C_i$ or $i+1\in C_i$. Therefore, it can be seen that
\begin{eqnarray*}
  \frac{\sum_{j\in N_i}c_{ij}}{\sum_{j\in C_i}c_{ij}}<-\frac{c_{ii}}{c_{ii-1}}=-\frac{c_{ii}}{c_{ii+1}}
  =-\frac{16\cos(\beta\pi)\Gamma(4-2\beta)+6\Gamma(3-\alpha)\tau^\alpha_nh^{-2\beta}(2^{4-2\beta}-
  8)}{4\cos(\beta\pi)\Gamma(4-2\beta)+6\Gamma(3-\alpha)\tau^\alpha_nh^{-2\beta}(3^{3-2\beta}-2^{5-2\beta}+7)},
\end{eqnarray*}
  indicating that $\sigma_2$ is independent of $h$ and $\tau_n$ by plugging \eqref{chp-03-015}, and thus prove the theorem.
\end{Proof}

\begin{Theorem}\label{thm-02}
  Let any C/F splitting be given. Under the condition \eqref{chp-03-015},
  there exist positive constants $\sigma_1$ and $\sigma_2$ independent of $h$ and $\tau_n$ and satisfying $\sigma_2>1>\sigma_1$,
  such that a uniform two-grid convergence is achieved as follows
\begin{eqnarray*}
  \|M_{h,H}\|_{1}\leq \sqrt{1-\sigma_1/\sigma_2}.
\end{eqnarray*}
\end{Theorem}

\begin{Proof}
  The proof of this result is straightforward and is based on Theorem A.4.1 and A.4.2 in \cite{s-001}, Lemma \ref{clp-03-01} and \ref{clp-03-02}.
\end{Proof}

Now observe from Theorem \ref{thm-02} that, despite the independence of $h$ and $\tau_n$,
$\sigma_2$ relies ruinously on $\theta$ in Step 1.1 of Algorithm \ref{alg-02}
due to the fact that $C^n_{h\tau}$ is nearly dense leading to a quite complicated adjacency graph.
In addition, it is found that $\sigma_2(\theta)$ may be much larger than 1 as $\theta$ approaches zero,
with that comes a sharp pullback in convergence rate.
Hence, an appropriate $\theta$ is a critical component of Algorithm \ref{alg-02} to handle fractional diffusion equations.

\begin{figure}[htp]
  \centerline{
  \includegraphics[scale=0.4]{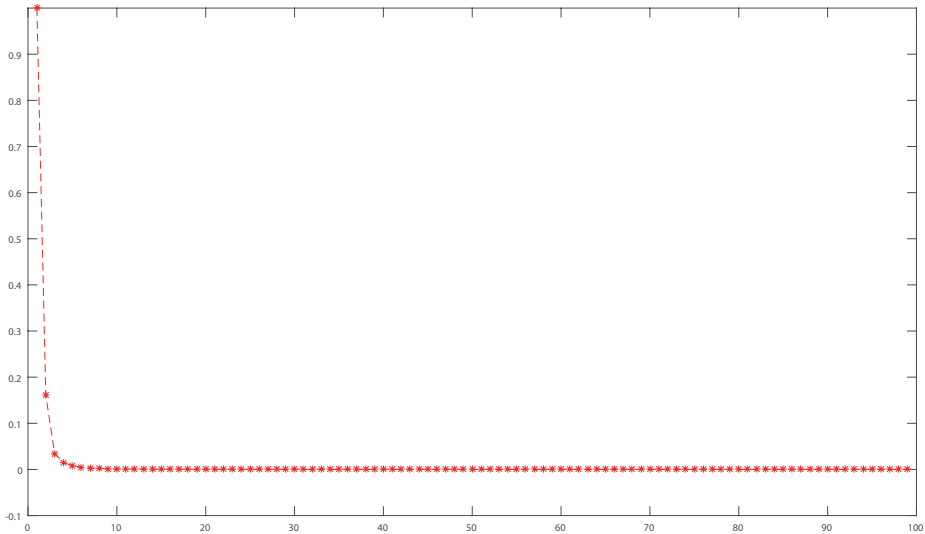}}
  \caption{Distribution of ratios $c_{1j}/c_{12}$, $j=2,3,\cdots,M-1$.}
  \label{fig-02}
\end{figure}

We now turn to reveal a reference formula on $\theta$. Note the heuristic that the distribution of ratios of off-diagonal elements
relative to the maximum absolute off-diagonal element (namely the minor diagonal element for $C^n_{h\tau}$) plays a major role
in the choice of $\theta$. Since $C^n_{h\tau}$ is a symmetric Toeplitz matrix from Corollary \ref{cor-01},
its first row involving all off-diagonal elements of $C^n_{h\tau}$ is deserved to be the representative row.
Taking $\beta=0.8$ as an example, Fig. \ref{fig-02} shows the distribution of the ratios $c_{1j}/c_{12}$ ($j\ge2$),
which reminds us of the attenuation in off-diagonal elements, states
\begin{eqnarray*}
 \frac{c_{13}}{c_{12}}\approx 0.160426,~~
 \frac{c_{1j}}{c_{12}} < \frac{c_{14}}{c_{12}} \approx 0.034394,~j=5,6,\cdots,M-1,
\end{eqnarray*}
and suggests that $c_{1j}$ ($j\ge4$) should be viewed as weak couplings (wouldn't be used for interpolation)
because they are less than 5\% of $c_{12}$. Besides, for a better complexity and higher efficiency,
only the nearest neighbors are potentially used to limit the interpolation matrix on each grid level
to at most 3 coefficients per row, although $c_{13}$ reaches around 16\% of $c_{12}$.
It thus appears that the strength-of-connection tolerance $\theta$ should be of the form
\begin{eqnarray}\label{chp-04-01}
  \theta = \frac{c_{13}}{c_{12}}+\epsilon_0,
\end{eqnarray}
where $\epsilon_0$ is some small number, which can be chosen to be $10^{-5}$ in one-dimensional realistic problems.

As is known, Algorithm \ref{alg-02} is much more expensive for well-conditioned problems than basic iterative techniques,
such as conjugate gradient (CG) or (plain) Jacobi iterative method. For the purpose of solving \eqref{chp-02-009}
in an optimal way, an adaptive AMG method is proposed below by combining Algorithm \ref{alg-02}, the reference formula \eqref{chp-04-01}
and the condition number estimation \eqref{chp-03-013}
in Corollary \ref{cor-03} as the clear distinction to adaptively pick an appropriate solver.

\begin{Algorithm}\label{alg-01} An adaptive AMG method $S_{ad}$ for the linear system \eqref{chp-02-009}.
  \begin{description}

    \item[Step 1] If the condition \eqref{chp-03-013} is unsatisfied, then goto Step 2,
    else set $S_{ad}$ as the CG or Jacobi iterative method;

    \item[Step 2] Set $S_{ad}$ as the classical AMG method described in Algorithm \ref{alg-02},
    with $\theta$ chosen via the reference formula \eqref{chp-04-01}.

  \end{description}
\end{Algorithm}

\section{Performance evaluation}

Let us illustrate the effectiveness of Algorithms \ref{alg-02} and \ref{alg-01}. Numerical experiments are performed in a 64 bit Fedora 18 platform,
double precision arithmetic on Intel Xeon (W5590) with 24.0 GB RAM, 3.33 GHz,
with an -O2 optimization parameter. In the following tables, dashed entries (-) indicate the solutions
either diverge or fail to converge after 1000 iterations,
\emph{Its} is the number of iterations until the stopping criterion $10^{-12}$ is reached,
$T_c$ represents the CPU time including both Setup and Solve phases with second as its unit,
$C_g$ and $C_o$ respectively denote grid and operator complexities, which are defined as
sums of the number of degrees of freedom and nonzero elements on all grid levels divided by those of the finest grid level,
and used as measures for memory requirements, aritmetic operations and the execution time in Setup and Solve phases.

\begin{Example}
  Comparisons of the classical AMG over CG and Jacobi iterative methods
  for the case when \eqref{chp-03-013} is satisfied with two different fractional orders.
\end{Example}

\begin{table}[htbp]
\footnotesize \centering\caption{Number of iterations and wall time for the case $\tau=h^2$.}\label{ctp-04-01}\vskip 0.1cm
\begin{tabular}{|c|c|c|c|c|c|c|c|c|c|c|c|c|}\hline
\multirow{3}{*}{$M$}&\multicolumn{6}{|c|}{$\alpha=\beta=0.6$}
&\multicolumn{6}{|c|}{$\alpha=\beta=0.8$} \\ \cline{2-13}
~&\multicolumn{2}{|c|}{Jacobi}&\multicolumn{2}{|c|}{CG}&\multicolumn{2}{|c|}{AMG}
 &\multicolumn{2}{|c|}{Jacobi}&\multicolumn{2}{|c|}{CG}&\multicolumn{2}{|c|}{AMG} \\ \cline{2-13}
~& \emph{Its} & $T_c$ & \emph{Its} & $T_c$ & \emph{Its} & $T_c$ & \emph{Its} & $T_c$ & \emph{Its} & $T_c$ & \emph{Its} & $T_c$ \\\hline

  32 & 18 & 1.78E-4 & 9  & 1.09E-4 & 4 & 5.22E-4 & 22 & 2.03E-4 & 11 & 1.27E-4 & 5 & 2.84E-4 \\ \hline
  64 & 18 & 3.90E-4 & 11 & 2.19E-4 & 4 & 7.79E-4 & 23 & 4.84E-4 & 13 & 1.99E-4 & 5 & 6.33E-4 \\ \hline
 128 & 19 & 1.31E-3 & 11 & 4.60E-4 & 4 & 2.52E-3 & 23 & 1.54E-3 & 13 & 5.26E-4 & 5 & 1.89E-3 \\ \hline
 256 & 19 & 4.69E-3 & 11 & 1.57E-3 & 4 & 9.73E-3 & 23 & 5.66E-3 & 13 & 1.82E-3 & 5 & 7.06E-3 \\ \hline
 512 & 19 & 2.61E-2 & 11 & 8.03E-3 & 4 & 5.49E-2 & 23 & 3.12E-2 & 13 & 9.56E-3 & 5 & 4.57E-2 \\ \hline
1024 & 19 & 1.95E-1 & 11 & 6.04E-2 & 4 & 1.73E-1 & 23 & 2.36E-1 & 12 & 6.53E-2 & 5 & 1.32E-1 \\ \hline
2048 & 19 & 3.98E-1 & 11 & 1.22E-1 & 4 & 9.39E-1 & 23 & 9.11E-1 & 12 & 1.32E-1 & 5 & 7.49E-1 \\ \hline
4096 & 19 &  3.03   & 11 & 9.25E-1 & 4 &  2.80   & 23 &  3.65   & 12 &  1.01   & 5 &  2.98   \\ \hline

\end{tabular}
\end{table}

As expected, the results in Table \ref{ctp-04-01} show that Jacobi, CG and AMG methods
are robust with respect to the mesh size and fractional order,
which indicates indirectly the correctness of \eqref{chp-02-007}.
In addition, CG method runs 3.28 and 3.03 times faster than Jacobi and AMG methods for $M=4096$ and $\alpha=\beta=0.6$, respectively.

\begin{Example}
  Comparisons between the classical AMG method and CG method for the case when \eqref{chp-03-013} is unsatisfied.
\end{Example}

\begin{table}[htbp]
\footnotesize
\centering\caption{Number of iterations and wall time for the case $\tau=1/32$.}\label{ctp-04-02}\vskip 0.1cm
\begin{tabular}{|c|c|c|c|c|c|c|c|c|c|c|c|c|c|}\hline
\multirow{3}{*}{$M$}&\multicolumn{4}{|c|}{$\beta=0.6$}&\multicolumn{4}{|c|}{$\beta=0.8$} &\multicolumn{4}{|c|}{$\beta=0.99$} \\ \cline{2-13}
~&\multicolumn{2}{|c|}{CG}&\multicolumn{2}{|c|}{AMG}&\multicolumn{2}{|c|}{CG}&\multicolumn{2}{|c|}{AMG}
&\multicolumn{2}{|c|}{CG}&\multicolumn{2}{|c|}{AMG} \\ \cline{2-13}
~& \emph{Its} & $T_c$ & \emph{Its} & $T_c$ & \emph{Its} & $T_c$ & \emph{Its} & $T_c$ & \emph{Its} & $T_c$ & \emph{Its} & $T_c$ \\\hline

 512 &  97 &  0.119 & 8 & 0.042 & 180 &  0.209 & 8 & 0.069 &     256   & 0.314 & 3 & 0.032 \\ \hline
1024 & 147 &  0.715 & 8 & 0.169 & 314 &  1.627 & 8 & 0.301 &     512   & 2.537 & 3 & 0.133 \\ \hline
2048 & 223 &  2.230 & 8 & 0.677 & 546 &  6.037 & 8 & 0.772 &  $>$1000  &   -   & 3 & 0.532 \\ \hline
4096 & 337 & 13.378 & 8 & 2.735 & 948 & 38.481 & 8 & 3.143 &  $>$1000  &   -   & 3 & 2.034 \\ \hline

\end{tabular}
\end{table}

As shown in Table \ref{ctp-04-02}, AMG method converges robustly regarding to the mesh size and may be weakly dependent of $\beta$,
while the number of iterations of CG method is quite unstable, and sometimes CG method even break down.
Furthermore AMG method runs 12.24 times faster than CG method for $M=4096$ and $\beta=0.8$.

\begin{table}[htbp]
\footnotesize
\centering\caption{Number of iterations and wall time for the case $\tau=h$.}\label{ctp-04-03}\vskip 0.1cm
\begin{tabular}{|c|c|c|c|c|c|c|c|c|}\hline
\multirow{3}{*}{$M$}&\multicolumn{4}{|c|}{$\alpha=0.2$, $\beta=0.6$} &\multicolumn{4}{|c|}{$\alpha=0.6$, $\beta=0.8$} \\ \cline{2-9}
~&\multicolumn{2}{|c|}{CG}&\multicolumn{2}{|c|}{AMG}&\multicolumn{2}{|c|}{CG}&\multicolumn{2}{|c|}{AMG} \\ \cline{2-9}
~& \emph{Its} & $T_c$ & \emph{Its} & $T_c$ & \emph{Its} & $T_c$ & \emph{Its} & $T_c$ \\\hline

 128 & 40  & 1.686E-3 & 8 & 2.868E-3 & 57  & 2.310E-3 & 7 & 2.873E-3 \\ \hline
 256 & 62  & 9.090E-3 & 8 & 1.571E-2 & 99  & 1.622E-2 & 7 & 1.081E-2 \\ \hline
 512 & 95  & 1.523E-1 & 8 & 6.976E-2 & 171 & 2.626E-1 & 7 & 4.808E-2 \\ \hline
1024 & 145 & 5.166E-1 & 8 & 2.916E-1 & 291 & 1.7895   & 7 & 2.063E-1 \\ \hline

\end{tabular}
\end{table}

Table \ref{ctp-04-03} shows the results of $\tau=h$.
Despite the advantage in computational cost and robustness over CG method,
AMG method is nearly independent of $\alpha$ and $\beta$ in this circumstance.
Meanwhile, by an investigation in terms of number of iterations in Tables \ref{ctp-04-02} and \ref{ctp-04-03},
CG method converges faster because of the improvement in condition number
from $\mathcal{O}(h^{-2\beta})$ to $\mathcal{O}(h^{\alpha-2\beta})$.

\begin{Example}
  Comparisons of $S_{ad}$ over the classical AMG and CG methods when the $i$-th time step size $\tau_i$ is chosen to be
\begin{eqnarray*}
  \tau_i = \left \{
    \begin{aligned}
    & h^2,~\quad i=1,\cdots,K_1 \\
    & 1/32,~i=K_1+1,\cdots,K_1+K_2
    \end{aligned}
  \right..
\end{eqnarray*}
\end{Example}

\begin{table}[htbp]
\footnotesize
\centering\caption{Comparisons among $S_{ad}$, CG and AMG.}\label{ctp-04-05}\vskip 0.1cm
\begin{tabular}{|c|c|c|c|c|c|c|c|c|c|c|c|c|}\hline
\multirow{3}{*}{$K_2$}&\multicolumn{6}{|c|}{$M=1024$, $K_1=K_2$}&\multicolumn{6}{|c|}{$M=2048$, $K_1=K_2$} \\ \cline{2-13}
~&\multicolumn{2}{|c|}{$S_\alpha$}&\multicolumn{2}{|c|}{CG}&\multicolumn{2}{|c|}{AMG}
&\multicolumn{2}{|c|}{$S_\alpha$}&\multicolumn{2}{|c|}{CG}&\multicolumn{2}{|c|}{AMG} \\ \cline{2-13}
~& \emph{Its} & $T_c$ & \emph{Its} & $T_c$ & \emph{Its} & $T_c$ & \emph{Its} & $T_c$ & \emph{Its} & $T_c$ & \emph{Its} & $T_c$ \\\hline

 25 &  459 & 2.09 &  7994 &  9.99 &  303 &  2.91 &  458 &  8.65 & 13681 &  71.76 &  304 & 12.03 \\ \hline
 50 &  909 & 4.11 & 15900 & 19.65 &  603 &  5.79 &  895 & 17.25 & 27196 & 141.81 &  604 & 23.79 \\ \hline
 75 & 1337 & 6.54 & 23748 & 30.74 &  903 &  8.60 & 1320 & 25.59 & 40637 & 212.48 &  904 & 35.76 \\ \hline
100 & 1760 & 8.54 & 31573 & 40.12 & 1201 & 11.92 & 1745 & 36.07 & 54038 & 292.03 & 1204 & 50.66 \\ \hline \hline

%
\multirow{3}{*}{$K_2$}&\multicolumn{6}{|c|}{$M=1024$, $K_1=3K_2$}&\multicolumn{6}{|c|}{$M=2048$, $K_1=3K_2$} \\ \cline{2-13}
~&\multicolumn{2}{|c|}{$S_{ad}$}&\multicolumn{2}{|c|}{CG}&\multicolumn{2}{|c|}{AMG}
&\multicolumn{2}{|c|}{$S_{ad}$}&\multicolumn{2}{|c|}{CG}&\multicolumn{2}{|c|}{AMG} \\ \cline{2-13}
~& \emph{Its} & $T_c$ & \emph{Its} & $T_c$ & \emph{Its} & $T_c$ & \emph{Its} & $T_c$ & \emph{Its} & $T_c$ & \emph{Its} & $T_c$ \\\hline

 25 &  987 &  2.86 &  8522 & 11.10 &  553 &  5.42 &  970 & 11.87 & 14193 &  79.22 &  554 & 21.86 \\ \hline
 50 & 1912 &  5.81 & 16903 & 22.57 & 1088 & 10.52 & 1895 & 22.66 & 28196 & 158.78 & 1091 & 47.96 \\ \hline
 75 & 2837 &  8.05 & 25248 & 31.73 & 1563 & 15.23 & 2820 & 33.83 & 42137 & 230.60 & 1566 & 66.73 \\ \hline
100 & 3760 & 10.82 & 33573 & 42.33 & 2036 & 20.03 & 3745 & 48.37 & 56038 & 294.82 & 2041 & 83.03 \\ \hline

\end{tabular}
\end{table}

We can observe from Table \ref{ctp-04-05} that $S_{ad}$ and AMG methods are fairly robust as to the mesh size,
roughly 10 and 6 on the average. Yet the average number of iterations of CG method varies from 85 to 142.
Moreover $S_{ad}$ has a considerable advantage over others in CPU time,
runs 1.72 and 6.09 times faster than AMG and CG methods for $M=2048$ and $K_2=100$.

\begin{Example}
  Analyze effects of the strength-of-connection tolerance $\theta$ on the performance of the classical AMG method.
\end{Example}

\begin{table}[htbp]
\footnotesize
\centering\caption{Effect of $\theta$ on the classical AMG when $M=512$.}\label{ctp-04-07}\vskip 0.1cm
\begin{tabular}{|c|c|c|c|c|c|c|c|c|c|c|}\hline
\multirow{2}{*}{$\theta$}&\multicolumn{4}{|c|}{$\beta=0.8$}&\multicolumn{4}{|c|}{$\beta=0.99$} \\ \cline{2-9}
~& \emph{Its} & $T_c$ & $C_g$ & $C_o$ & \emph{Its} & $T_c$ & $C_g$ & $C_o$ \\\hline

 0.0001  & 293 & 1.952    & 1.037 & 1.001 & 103 & 6.784E-1 & 1.170 & 1.021 \\ \hline
 0.001   & 83  & 5.618E-1 & 1.098 & 1.008 & 60  & 4.189E-1 & 1.498 & 1.124 \\ \hline
0.00684  & 31  & 1.797E-1 & 1.202 & 1.029 & 32  & 1.374E-1 & 1.652 & 1.147 \\ \hline
0.00685  & 31  & 1.789E-1 & 1.202 & 1.029 & 3   & 3.255E-1 & 1.975 & 1.331 \\ \hline
 0.01    & 23  & 1.020E-1 & 1.247 & 1.041 & 3   & 3.301E-2 & 1.975 & 1.331 \\ \hline
 0.1     & 13  & 6.209E-2 & 1.489 & 1.124 & 3   & 4.356E-2 & 1.975 & 1.331 \\ \hline
 0.16042 & 13  & 6.037E-2 & 1.489 & 1.124 & 3   & 4.118E-2 & 1.975 & 1.331 \\ \hline
 0.16043 & 7   & 4.699E-2 & 1.975 & 1.331 & 3   & 3.158E-2 & 1.975 & 1.331 \\ \hline
 0.25    & 7   & 4.736E-2 & 1.975 & 1.331 & 3   & 4.353E-2 & 1.975 & 1.331 \\ \hline
 0.5     & 7   & 4.710E-2 & 1.975 & 1.331 & 3   & 4.354E-2 & 1.975 & 1.331 \\ \hline

\end{tabular}
\end{table}

\begin{table}[htbp]
\footnotesize
\centering\caption{Effect of $\theta$ on the classical AMG when $M=2048$.}\label{ctp-04-06}\vskip 0.1cm
\begin{tabular}{|c|c|c|c|c|c|c|c|c|c|c|}\hline
\multirow{2}{*}{$\theta$}&\multicolumn{4}{|c|}{$\beta=0.8$}&\multicolumn{4}{|c|}{$\beta=0.99$} \\ \cline{2-9}
~& \emph{Its} & $T_c$ & $C_g$ & $C_o$ & \emph{Its} & $T_c$ & $C_g$ & $C_o$ \\\hline

 0.0001 & 335 & 23.621 & 1.038 & 1.001 & 123 &   6.310  & 1.170 & 1.021 \\ \hline
 0.001  & 107 &  9.024 & 1.100 & 1.008 & 102 &   5.597  & 1.497 & 1.125 \\ \hline
0.00684 & 33  &  3.096 & 1.207 & 1.029 & 33  &   2.055  & 1.662 & 1.148 \\ \hline
0.00685 & 33  &  3.129 & 1.207 & 1.029 & 3   & 7.314E-1 & 1.993 & 1.333 \\ \hline
 0.01   & 26  &  2.551 & 1.249 & 1.042 & 3   & 5.369E-1 & 1.993 & 1.333 \\ \hline
 0.1    & 15  &  1.691 & 1.497 & 1.125 & 3   & 5.448E-1 & 1.993 & 1.333 \\ \hline
 0.1603 & 15  &  1.690 & 1.497 & 1.125 & 3   & 5.372E-1 & 1.993 & 1.333 \\ \hline
 0.1604 & 8   &  1.217 & 1.993 & 1.333 & 3   & 5.294E-1 & 1.993 & 1.333 \\ \hline
 0.2    & 8   &  1.219 & 1.993 & 1.333 & 3   & 5.294E-1 & 1.993 & 1.333 \\ \hline
 0.25   & 8   &  1.217 & 1.993 & 1.333 & 3   & 5.299E-1 & 1.993 & 1.333 \\ \hline

\end{tabular}
\end{table}

It is seen from Tables \ref{ctp-04-07} and \ref{ctp-04-06} that there is a unique threshold $\theta_0$ independent of $h$
which guarantees the robustness of the classical AMG method, and makes number of iterations of the classical AMG monotonically decreasing
when $\theta<\theta_0$, or even the classical AMG possibly diverge when $\theta$ is small enough,
e.g., $\theta_0=0.16043$ and $\theta_0=0.00685$ for cases $\beta=0.8$ and $\beta=0.99$.
By direct calculations, we have $c_{13}/c_{12}\approx 0.160426$ and $c_{13}/c_{12}\approx 0.006846$.
Utilizing the relation \eqref{chp-04-01} and $\epsilon_0=10^{-5}$,
the corresponding values of $\theta$ are respectively larger than those of $\theta_0$.
This confirms the reasonability of the reference formula \eqref{chp-04-01}.

\section{Conclusion}

In this paper, we propose the variational formulation for a class of time-space Caputo-Riesz fractional diffusion equations,
prove that the resulting matrix is a symmetric Toeplitz matrix, an M-matrix by appending a very weak constraint
and its condition number is bounded by $1+\mathcal{O}(\tau^\alpha h^{-2\beta})$,
introduce the classical AMG method and prove rigorously that its convergence rate is independent of time and space step sizes,
provide explicitly a reference formula of the strength-of-connection tolerance to guarantee the robustness and predictable behavior of AMG method in all cases,
and develop an adaptive AMG method via our condition number estimation to decrease the computation cost.
Numerical results are all in conformity with the theoretical results, and verify the reasonability of the reference formula
and the considerable advantage of the proposed adaptive AMG algorithm over other traditional iterative methods, e.g. Jacobi,
CG and the classical AMG methods.

\section*{Acknowledgments}

This work is under auspices of National Natural Science Foundation of China (11571293, 11601460, 11601462)
and the General Project of Hunan Provincial Education Department of China (16C1540, 17C1527).

\end{document}